\input amstex
\magnification=\magstep1
\baselineskip=13pt
\documentstyle{amsppt}
\vsize=8.7truein
\CenteredTagsOnSplits
\NoRunningHeads
\def\today{\ifcase\month\or
  January\or February\or March\or April\or May\or June\or
  July\or August\or September\or October\or November\or December\fi
  \space\number\day, \number\year}
\def\Pr{\bold{P\ }}
\def\EE{\bold{E\ }}
\def\AA{\Cal{A}}
\def\per{\operatorname{per}}
\def\sgn{\operatorname{sgn}}
\def\Mat{\operatorname{Mat}}
\def\tr{\operatorname{Tr}}
\def\Cdet{\operatorname{Cdet}}
\def\sdet{\operatorname{sdet}}
\NoRunningHeads
\topmatter
\title New Permanent Estimators via Non-Commutative Determinants \endtitle
\author Alexander Barvinok \endauthor
\address Department of Mathematics, University of Michigan, Ann Arbor,
MI 48109-1109 \endaddress
\email barvinok$\@$math.lsa.umich.edu  \endemail
\date  July 2000 \enddate
\thanks This research  was partially supported by NSF Grant DMS 9734138.
\endthanks
\abstract We introduce a new notion of the determinant,
called symmetrized determinant, for a square matrix 
with the entries in an associative algebra $\AA$. The monomial expansion
of the symmetrized determinant is obtained from the standard 
expansion of the commutative determinant by averaging the products of entries 
of the matrix in all possible orders. We show that for any fixed 
finite-dimensional associative algebra $\AA$, the symmetrized 
determinant of an $n\times n$ matrix with the entries in $\AA$ can be 
computed in polynomial in $n$ time (the degree of the polynomial is 
linear in the dimension of $\AA$). 
Then, for every associative 
algebra $\AA$ endowed with a scalar product and unbiased probability 
measure, we construct a randomized polynomial time algorithm to 
estimate the permanent of non-negative matrices. We conjecture that if 
$\AA=\Mat(d, {\Bbb R})$ is the algebra of $d\times d$ real matrices
endowed with the standard scalar product and Gaussian measure,
the algorithm approximates the permanent of a non-negative $n \times n$ 
matrix within $O(\gamma_d^n)$ factor, where 
$\lim_{d \longrightarrow +\infty} \gamma_d=1$.
Finally, we provide some informal arguments 
why the conjecture might be true.    
\endabstract
\endtopmatter
\document

\head 1. Introduction \endhead

\subhead (1.1) Permanent and determinant \endsubhead
Let $A=(a_{ij})$ be an $n \times n$ real matrix and let 
$S_n$ denote the symmetric group of all permutations $\sigma$ 
of the set $\{1, \ldots, n\}$. The number 
$$\per A=\sum_{\sigma \in S_n} \prod_{i=1}^n a_{i \sigma(i)}$$
is called the {\it permanent} of $A$. We will be interested in the 
case when $A$ is non-negative: $a_{ij} \geq 0$ for $i,j=1, \ldots, n$. 

To compute or to approximate 
the permanent efficiently is one of the most intriguing problems,
see, for example, [Jerrum and Sinclair 89], [Karmarkar {\it et al.} 93],
Chapter 18 of [Papadimitriou 94], [Linial {\it et al.} 98],
[Barvinok 97] and [Barvinok 99]. 

On the other hand, the determinant
$$\det A=\sum_{\sigma \in S_n} (\sgn \sigma) 
\prod_{i=1}^n a_{i \sigma(i)}$$
can be computed using $O(n^3)$ arithmetic operations (even when 
$A$ is a matrix over a commutative ring), see, for example, Sections 11.1
and 15.1 of 
[Papadimitriou 94].

Hence it seems natural to try to use determinants to approximate 
the permanent.
\subhead (1.2) Permanent estimators \endsubhead
The following approach is due to Godsil and Gutman and 
independently to Girko (see, for example, 
Chapter 8 of  [Lov\'asz and Plummer 86] and Section 3 of Chapter 2 of
[Girko 90].
Let us choose a probability distribution $\mu$ on 
${\Bbb R}$ with the properties that 
$$\EE x = \int_{\Bbb R} x \ d \mu(x) =0 \quad  \text{and} \quad  
\EE x^2 =\int_{\Bbb R} x^2 \ d \mu(x)=1.$$ Let us sample $n^2$ numbers 
$u_{ij}: i,j=1, \ldots, n$ independently and at random from $\mu$.
Let $B=(b_{ij})$ be the matrix defined by
$$b_{ij}= u_{ij}\sqrt{a_{ij}}$$ and 
let 
$$\alpha =(\det B)^2.$$
Then $\alpha$ is a random variable and it turns 
out that $\EE \alpha =\per A$. Hence, in principle, $\per A$ can be 
approximated by averaging sufficiently many determinants.
It turns out that by extending the ground field, one can 
ensure a better concentration of $\alpha$ around its expectation
and hence a better computational complexity of the approximation.

In  [Karmarkar {\it et al.} 93]
it was shown that if one allows $u_{ij}$ to be complex 
numbers (we must define then $\alpha=|\det B|^2$), one can make 
the variance of $\alpha$ smaller. More precisely, the variance 
of $\alpha$ for the distribution $\mu$ that chooses the cubic roots 
of unity with the probability $1/3$ each is exponentially 
smaller in the worst case
than that for $\mu$ that chooses $-1$ and $1$ with the 
probability $1/2$ each. In [Barvinok 99], it was shown that 
if $\mu$ is a Gaussian distribution in ${\Bbb R}$ then with 
high probability a random value of $\alpha$ approximates its 
expectation within a $c^n$ factor for $c \approx 0.28$; similar 
behavior is observed when $\mu$ is a complex Gaussian distribution,
but in the complex case the constant $c$ gets 
better: $c \approx 0.56$. Moreover, if $\mu$ is a quaternionic Gaussian 
distribution (in which case $\alpha$ should be the  
{\it Study determinant} of $B$, that is the 
determinant of the complexification of $B$, 
 see, for example, [Aslaksen 96]), the 
constant $c$ gets even better: $c \approx 0.76$.

This suggests that it may be of interest to construct more 
determinant-type estimators, as they could provide still better 
concentration properties of $\alpha$.  

Let us summarize some useful properties of the estimator.
The estimator $\alpha=\alpha(u_{ij}): i,j=1, \ldots, n$
is a function of $n^2$ independent and identically distributed 
random variables $u_{ij}$, such that:
\bigskip
\noindent (1.2.1) For each choice of $u_{ij}$, the value 
of $\alpha(u_{ij})$ can be computed 
using a polynomial in $n$ number of arithmetic operations;
\medskip
\noindent (1.2.2) The expected value of $\alpha$ is $\per A$;
\medskip
\noindent (1.2.3) For any choice of $(u_{ij})$, the value of 
$\alpha(u_{ij})$ is a non-negative real number;
\medskip
\noindent(1.2.4) Let us fix all $u_{ij}$ except those in one row $i_0$ 
(resp. in one column $j_0$). Then $\alpha$ is a quadratic form 
in $\{u_{i_0j}: j=1, \ldots, n\}$  
(resp. in $\{ u_{i j_0}: i=1, \ldots, n \}$).
\bigskip
In this paper, for every associative finite-dimensional algebra 
$\AA$ we construct an estimator with the properties (1.2.1)--(1.2.4).
\subhead (1.3) Non-commutative determinants \endsubhead
Exploring analogues of determinants over non-commutative rings and 
algebras is a very old topic, see, for example, [Aslaksen 96] for a survey 
and [Gelfand and Retakh 97] for new developments. Interestingly, most important
non-commutative determinants, such as the Dieudonn\'e determinant, 
quasideterminants of Gelfand and Retakh
and the Moore determinant of a 
Hermitian quaternionic matrix turn out to be computationally efficient.
However, with the notable exception of the Moore determinant,
they are not close enough to the monomial expansion of the commutative  
determinant to produce a permanent estimator. On the other hand, 
the most straightforward version due to Cayley
$$\Cdet A=\sum_{\sigma \in S_n} (\sgn \sigma) a_{1 \sigma(1)} \cdots 
a_{n \sigma(n)}$$
for an $n \times n$ matrix $A=(a_{ij})$ with the entries $a_{ij}$ in 
an associative algebra $\AA$ would have suited our purpose
perfectly well had there  
been any reason to believe that this expression can be 
efficiently (in polynomial in $n$ time) computed, see also Section 5.2. 

The ambiguity that prevents the straightforward extension of the 
determinant to non-commutative algebras is that there is no natural 
choice for the order of the factors in the product 
$\prod_{i=1}^n a_{i \sigma(i)}$. We resolve this ambiguity by taking
the average of the products in all possible $n!$ orders. Hence for 
an $n \times n$ matrix $A=(a_{ij})$ with the entries in an 
associative algebra $\AA$ over a field ${\Bbb F}$ of characteristic 0, 
we write
$$\sdet A={1 \over n!} \sum_{(\sigma, \tau) \in S_n \times S_n} 
(\sgn \sigma)(\sgn \tau) a_{\sigma(1) \tau(1)} \cdots a_{\sigma(n) \tau(n)}$$
(``sdet'' stands for the ``symmetrized determinant'').
It is easy to see that if $\AA$ is commutative, we get the standard 
determinant. We prove that for any fixed finite-dimensional algebra 
$\AA$, the value of $\sdet A$ for an $n\times n$ matrix 
over $\AA$ can be computed in a polynomial in $n$ time. More precisely,
if the dimension of $\AA$ as an ${\Bbb F}$-vector space is $r$, 
$\sdet A$ can be computed in $O(n^{r+3})$ time.   
    
The paper is organized as follows.
In Section 2, we discuss mixed discriminants, which are crucial for 
our proof in Section 3 of the polynomial time computability of the 
symmetrized determinant. 
In Section 4, for every finite-dimensional associative algebra
$\AA$ over ${\Bbb R}$ endowed with a 
scalar product and an unbiased probability distribution $\mu$, we construct 
a permanent estimator which satisfies (1.2.1)--(1.2.4). In Section 4, 
we conjecture that if $\AA=\Mat(d, {\Bbb R})$ is the algebra of 
$d \times d$ matrices endowed with the standard scalar product and 
Gaussian probability measure, then the algorithm approximates the permanent of 
an $n \times n$ matrix within a $O(\gamma_d^n)$ factor, where
$\lim_{d \longrightarrow +\infty} \gamma_d=1$. We also 
provide some intuitive argument supporting the conjecture. 
 
\head 2. Preliminaries: Mixed Discriminants \endhead

\definition{(2.1) Definition} Let $A_1, \ldots, A_n$ be  
$n \times n$ matrices over a field ${\Bbb F}$ of 
characteristic 0.
We write $A_k=(a_{ij}^k)$, where $a_{ij}^k \in {\Bbb F}$
for $i,j=1, \ldots, n$ and $k=1, \ldots, n$ ($k$ is the index,
not the power). Let $t_1, \ldots, t_n \in {\Bbb F}$ be variables. The 
expression $\det\bigl(t_1 A_1 + \ldots + t_n A_n \bigr)$ is a homogeneous 
polynomial in $t_1, \ldots, t_n$ of degree $n$ and its
normalized coefficient
$$D(A_1, \ldots, A_n)={1 \over n!} {\partial^n \over 
\partial t_1 \ldots \partial t_n} \det\bigl(t_1 A_1 + \ldots + t_n A_n\bigr)$$
is called the {\it mixed discriminant} of $A_1, \ldots, A_n$.
In terms of the entries $(a_{ij}^k)$ of the matrices $A_1, \ldots, A_n$, 
the mixed discriminant can be written as 
$$D(A_1, \ldots, A_n)={1 \over n!} \sum_{(\sigma, \tau) \in S_n \times S_n} 
(\sgn \sigma) (\sgn \tau) \prod_{k=1}^n a^k_{\sigma(k) \tau(k)}. 
\tag2.1.1$$  
\enddefinition
Mixed discriminants are symmetric, that is,
$$D(A_1, \ldots, A_n)=D(A_{\phi(1)}, \ldots, A_{\phi(n)}) \tag2.2$$ 
for any permutation $\phi: \{1, \ldots, n\} \longrightarrow 
\{1, \ldots, n\}$. For various properties of mixed discriminants,
see, for example, Section 5.2 of [Bapat and Raghavan 97].

We are particularly interested in the situation when the number of 
different
matrices among $A_1, \ldots, A_n$ is small. The following result is 
Lemma 9.3 from [Barvinok 97]. For the sake of completeness, we
present its proof here.
\proclaim{(2.3) Lemma} Let $k_1, \ldots, k_r$ be non-negative integers
such that $k_1 + \ldots + k_r=n$ and let $A_1, \ldots, A_r$ be $n \times n$ 
matrices. Then 
$$\split &D\bigl(\underbrace{A_1, \ldots, A_1}_{\text{$k_1$ copies}}, \ldots, 
\underbrace{A_r, \ldots, A_r}_{\text{$k_r$ copies}} \bigr) \\ =
&{(-1)^n \over n!}
\sum\Sb 0 \leq m_1 \leq k_1 \\ \ldots \\ 0 \leq m_r \leq k_r 
\endSb (-1)^{m_1 + \ldots +m_r}
{k_1 \choose m_1} \cdots {k_r \choose m_r} \det \bigl(m_1 A_1 + \ldots 
+m_r A_r \bigr). \endsplit $$
\endproclaim
\demo{Proof} For a subset $\omega \subset \{1, \ldots, n\}$, let 
$$t_i(\omega)=\cases 1 &\text{if\ } i \in \omega \\ 0 &\text{if \ }
i \notin \omega \endcases$$
and let ${\bold t}_{\omega}=\bigl(t_1(\omega), \ldots, t_n(\omega)\bigr)$
be the indicator of $\omega$.
One can observe that if $p$ is a homogeneous polynomial of 
degree $n>0$ in $n$ variables $t_1, \ldots, t_n$, then 
$${ \partial^n \over \partial t_1 \cdots \partial t_n} p= 
(-1)^n \sum_{\omega \subset \{1, \ldots, n \}} (-1)^{|\omega|} 
p({\bold t}_{\omega}).$$
Indeed, it suffices to check the identity for monomials 
$t_1^{\alpha_1} \cdots t_n^{\alpha_n}$. If some $\alpha_i=0$ then 
the right hand side is 0 since the terms corresponding to 
$\omega \cup \{i\}$ and 
$\omega \setminus \{i\}$ annihilate each other. For the monomial 
$t_1 \cdots t_n$ the right hand side is 1 since the only 
non-zero term is for $\omega=\{1, \ldots, n\}$. 

Given $n \times n$ matrices $A_1, \ldots, A_n$, let us apply 
the above identity to 
$$p(t_1, \ldots, t_n)=\det(t_1 A_1 + \ldots + t_n A_n).$$
Suppose that the 
 set $\{A_1, \ldots, A_n\}$ consists of $k_1$ copies of $A_1$, $\ldots$,
$k_r$ copies of $A_r$, where 
$A_1, \ldots, A_r$ are distinct, 
and let $S_i=\bigl\{j: A_j$ is a copy of $A_i \bigr\}$
for $i=1, \ldots, r$. Then $|S_i|=k_i$,
 $S_1 \cup  \ldots \cup S_r$ is a partition of the set 
$\{1, \ldots, n\}$ and  
$$p({\bold t}_{\omega})=\det(m_1 A_1 + \ldots + m_r A_r),
\quad \text{where} \quad m_i=|\omega \cap S_i|.$$
Moreover, for given $m_1, \ldots, m_r$, there are exactly
$\displaystyle{k_1 \choose m_1} \cdots {k_r \choose m_r}$ subsets 
$\omega$ with 
$m_i=|\omega \cap S_i|$, since for each $i=1, \ldots, r$,
we have to choose $m_i$ elements of $\omega$ from
$S_i$ for all $i=1, \ldots, r$ independently. The proof now follows.
{\hfill \hfill \hfill} \qed
\enddemo
\proclaim{(2.4) Corollary} Suppose that $r$ is fixed. Given $n \times n$ 
matrices $A_1, \ldots, A_r$, computing 
$$D\bigl(\underbrace{A_1, \ldots, A_1}_{\text{$k_1$ copies}}, \ldots, 
\underbrace{A_r, \ldots, A_r}_{\text{$k_r$ copies}} \bigr)$$
using the formula of Lemma 2.3, takes $O(n^{r+3})$ arithmetic operations.
\endproclaim
\demo{Proof} The number of summands does not exceed $(n+1)^r$ (a better 
estimate is $(n/r+1)^r$) and 
the determinant of an $n \times n$ matrix can be computed using 
$O(n^3)$ arithmetic operations.
{\hfill \hfill \hfill} \qed
\enddemo                                                               

\head 3. Symmetrized Determinant of a Matrix over an Algebra \endhead 

Let $\AA$ be an associative algebra over ${\Bbb F}$, 
where ${\Bbb F}$ is a field of characteristic 0. Hence $\AA$ is 
a vector space over ${\Bbb F}$ with an addition ``+'' and associative
(but not necessarily commutative) multiplication $\cdot$.
For example, one can choose $\AA$ to be the algebra $\Mat(d, {\Bbb F})$ of all 
$d \times d$ matrices over ${\Bbb F}$. 
We will assume that as a vector space, $\AA$ is finite-dimensional,
$\dim_{\Bbb F} \AA < \infty$. 
\definition{(3.1) Definition} Let $A=(a_{ij})$,
$a_{ij} \in \AA$ be an $n \times n$ matrix over $\AA$. 
We call 
$$\sdet A={1 \over n!} 
\sum_{(\sigma, \tau) \in S_n \times S_n} 
(\sgn \sigma) (\sgn \tau) a_{\sigma(1) \tau(1)} a_{\sigma(2) \tau(2)}
\cdots a_{\sigma(n) \tau(n)}$$
the {\it symmetrized determinant} of $A$. In other words, $\sdet A$ 
is obtained by taking a diagonal of $A$ (that is, picking one entry 
from each row and column of $A$), multiplying the entries on the 
diagonal in all possible orders, taking the average of the resulting 
$n!$ products and adding that average with the appropriate sign found
by the same rule as for the usual commutative determinant. 
Hence $\sdet A \in \AA$. 
We denote by $\Mat(n, \AA)$ the set of $n \times n$ matrices $A=(a_{ij})$ with
$a_{ij} \in \AA$.
\enddefinition
\remark{Remark} If $\AA$ is commutative, the expressions 
$${1 \over n!} \sum_{(\sigma, \tau) \in S_n \times S_n} 
(\sgn \sigma) (\sgn \tau) a_{\sigma(1) \tau(1)} 
\cdots a_{\sigma(n) \tau(n)}$$ 
and 
$$\sum_{\sigma \in S_n} (\sgn \sigma) a_{1 \sigma(1)} \cdots 
a_{n \sigma(n)}$$
coincide but if $\AA$ is not commutative, they may differ.
\endremark
\proclaim{(3.2) Theorem} Let $e_1, \ldots, e_r$ span $\AA$ as 
an ${\Bbb F}$-vector space. Let $A=(a_{ij})$ be an $n \times n$ matrix 
over $\AA$, so 
$$A=A_1 e_1 + \ldots + A_r e_r,$$
where $A_1, \ldots, A_r$ are $n \times n$ matrices with the entries in 
${\Bbb F}$. 

For an $r$-tuple $k_1, \ldots, k_r$ of non-negative 
integers such that $k_1 + \ldots + k_r=n$,
let $\Phi(k_1, \ldots, k_r)$ be the set of all 
maps $\phi: \{1, \ldots, n\} \longrightarrow \{1, \ldots, r\}$ such
that $|\phi^{-1}(1)|=k_1, \ldots, |\phi^{-1}(r)|=k_r$.
Let us define elements $u(k_1, \ldots, k_r) \in \AA$ by 
$$u(k_1, \ldots, k_r)=\sum_{\phi \in \Phi(k_1, \ldots, k_r)} 
e_{\phi(1)} e_{\phi(2)} \cdots e_{\phi(n)}.$$
Then 
$$\sdet A=\sum \Sb k_1, \ldots, k_r \geq 0 \\ 
k_1 + \ldots + k_r=n \endSb 
D\bigl(\underbrace{A_1, \ldots, A_1}_{\text{$k_1$ copies}}, \ldots, 
\underbrace{A_r, \ldots, A_r}_{\text{$k_r$ copies}} \bigr)
u(k_1, \ldots, k_r).$$ 
\endproclaim
\demo{Proof}  Let $A_k=(a_{ij}^k)$, where $a_{ij}^k \in {\Bbb F}$
for $i,j=1,\ldots, n$ and $k=1, \ldots, r$. 
Hence $a_{ij} = \sum_{k=1}^r a_{ij}^k e_k$.

By Definition 3.1, 
$$\split \sdet A &= {1 \over n!} \sum_{(\sigma, \tau) \in S_n \times S_n} 
(\sgn \sigma) (\sgn \tau) a_{\sigma(1) \tau(1)} a_{\sigma(2) \tau(2)}
\cdots a_{\sigma(n) \tau(n)} \\ & =
 {1 \over n!} \sum_{(\sigma, \tau) \in S_n \times S_n} 
(\sgn \sigma) (\sgn \tau) \Bigl( \sum_{k=1}^r a_{\sigma(1) \tau(1)}^k e_k
 \Bigr) \ldots \Bigl( \sum_{k=1}^r a_{\sigma(n) \tau(n)}^k e_k \Bigr)
\endsplit$$
Let $\Phi$ be the set of all maps 
$\phi:\{1, \ldots, n\} \longrightarrow \{1, \ldots, r\}$.
Then 
$$\Bigl(\sum_{k=1}^r a_{\sigma(1) \tau(1)}^k 
e_k \Bigr) \cdots \Bigl( \sum_{k=1}^r a_{\sigma(n) \tau(n)}^k e_k \Bigr)=
\sum_{\phi \in \Phi} a_{\sigma(1) \tau(1)}^{\phi(1)} \cdots 
a_{\sigma(n) \tau(n)}^{\phi(n)} e_{\phi(1)} \cdots e_{\phi(n)}.$$
Hence
$$\split \sdet A &= {1 \over n!} \sum_{(\sigma, \tau) \in S_n \times S_n} 
(\sgn \sigma) (\sgn \tau) \sum_{\phi \in \Phi} 
a_{\sigma(1) \tau(1)}^{\phi(1)} \cdots 
a_{\sigma(n) \tau(n)}^{\phi(n)} e_{\phi(1)} \cdots e_{\phi(n)} \\ &=
\sum_{\phi \in \Phi} {1 \over n!} 
\sum_{(\sigma, \tau) \in S_n \times S_n} 
(\sgn \sigma) (\sgn \tau) a_{\sigma(1)\tau(1)}^{\phi(1)} \cdots 
 a_{\sigma(n) \tau(n)}^{\phi(n)} 
 e_{\phi(1)} \cdots e_{\phi(n)} \\ 
&= \sum_{\phi \in \Phi} 
D\bigl(A_{\phi(1)}, \ldots, A_{\phi(n)}\bigr)
e_{\phi(1)} \cdots e_{\phi(n)}. \endsplit $$ 
Now, set $\Phi$ is a disjoint union of the sets $\Phi(k_1, \ldots, k_r)$ and 
by the symmetry of the mixed discriminant 
(see (2.2)), for any $\phi \in \Phi(k_1, \ldots, k_r)$ we have 
$$D(A_{\phi(1)}, \ldots A_{\phi(n)})=
D\bigl(\underbrace{A_1, \ldots, A_1}_{\text{$k_1$ copies}}, \ldots, 
\underbrace{A_r, \ldots, A_r}_{\text{$k_r$ copies}} \bigr).$$
Summarizing, we get
$$\sdet A=\sum \Sb k_1, \ldots, k_r \geq 0 \\ k_1 + \ldots 
+ k_r =n \endSb D\bigl(\underbrace{A_1, \ldots, A_1}_{\text{$k_1$ copies}},
 \ldots, 
\underbrace{A_r, \ldots, A_r}_{\text{$k_r$ copies}} \bigr)
\sum_{\phi \in \Phi(k_1, \ldots, k_r)} e_{\phi(1)} \cdots e_{\phi(n)}$$
and the proof follows. 
{\hfill \hfill \hfill} \qed
\enddemo
Next, we present an algorithm for computing the 
symmetrized determinant of a matrix over $\AA$. We assume that there is 
a basis $e_1, \ldots, e_r$ of $\AA$ as a vector space 
and that the entries $a_{ij}$ of $A \in \Mat(n, \AA)$ 
are given by their coefficients 
in the basis: $a_{ij}=a_{ij}^1 e_1 + \ldots + a_{ij}^r e_r$. 
The multiplication in $\AA$ is described by the structural 
constants $\{\beta_{ij}^k\}$ that are the scalars from ${\Bbb F}$ such that
$$e_i e_j =\sum_{k=1}^r \beta_{ij}^k e_k. \leqno(3.3)$$
We need a simple Lemma.
\proclaim{(3.4) Lemma} Let $u(k_1, \ldots, k_r) \in \AA$ be the elements 
defined in Theorem 3.2. Let $I=\{i: k_i > 0 \}$.
Then 
$$u(k_1, \ldots, k_r)=\sum_{i \in I} u(k_1, \ldots, k_{i-1}, k_i-1,
k_{i+1}, \ldots, k_r) e_i.$$
\endproclaim
\demo{Proof} For $\phi \in \Phi(k_1, \ldots, k_r)$ we have
$$e_{\phi(1)} \cdots e_{\phi(n)}=\Bigl(e_{\phi(1)} \cdots e_{\phi(n-1)} 
\Bigr) e_{\phi(n)},\qquad \text{where} \quad \phi(n) \in I.$$
If $\phi(n)=i$ then the restriction
$\phi': \{1, \ldots, n-1\} \longrightarrow \{1, \ldots, r\}$ belongs to 
the set $\Phi(k_1, \ldots, k_{i-1}, k_i-1, k_{i+1}, \ldots, k_r)$.
Vice versa, every
$\phi' \in \Phi(k_1, \ldots, k_{i-1}, k_i-1, k_{i+1}, \ldots, k_r)$
extends to $\phi \in \Phi(k_1, \ldots, k_r)$ by letting 
$\phi(n)=i$. The proof now follows.
{\hfill \hfill \hfill} \qed
\enddemo
\specialhead(3.5) Algorithm for computing $\sdet A$ for $A \in \Mat(n, \AA)$
\endspecialhead
\noindent{\bf Input:} An algebra $\AA$ given by the structural constants 
$\{\beta_{ij}^k \}$ in some ${\Bbb F}$-basis $e_1, \ldots, e_r$ of 
$\AA$ and an $n \times n$ matrix $A \in \Mat(n, \AA)$ given by 
$n \times n$ matrices $A_1, \ldots, A_r \in \Mat(n, {\Bbb F})$ 
such that 
$A=A_1 e_1 + \ldots + A_r e_r$.
\medskip
\noindent{\bf Output:} The element $\sdet A \in \AA$ given by 
its coefficients in the basis $e_1, \ldots, e_r$. 
\medskip
\noindent{\bf Algorithm:} First, for all $r$-tuples of non-negative 
integers $k_1, \ldots, k_r$ such that $k_1 + \ldots + k_r=n$, using 
the identity of Lemma 2.3, we compute the mixed discriminant 
$$D\bigl(\underbrace{A_1, \ldots, A_1}_{\text{$k_1$ copies}},
 \ldots, 
\underbrace{A_r, \ldots, A_r}_{\text{$k_r$ copies}} \bigr).$$
Now, using Lemma 3.4 recursively, we compute $u(k_1, \ldots, k_r)$.
We start with 
$$u(0, \ldots, 0, \underset \text{$i$-th position} \to {1,} 
0, \ldots, 0)=e_i$$ 
and applying Lemma 3.4 and (3.3), successively compute the 
expansions of \break $u(k_1, \ldots, k_r)$ with 
$k_1 + \ldots +k_r=1, \ldots, n$ in the basis $e_1, \ldots, e_r$.
Finally, we use Theorem 3.2 to compute $\sdet A$. 
\bigskip 
An important observation is that if $\dim_{\Bbb F} \AA=r$ is fixed, 
the algorithm has polynomial time complexity.
\proclaim{(3.6) Theorem} For a fixed $r$, given an $n \times n$ 
matrix $A \in \Mat(\AA, n)$, Algorithm 3.5 computes $\sdet A \in \AA$ using 
$O(n^{r+3})$ arithmetic operations. 
\endproclaim
\demo{Proof} Theorem 3.2 and Lemma 3.4 imply that the algorithm indeed
returns the correct value. The number of all $r$-tuples 
of non-negative integers $k_1, \ldots, k_r$ such that 
$k_1 + \ldots + k_r \leq n$ is ${n+r \choose r}$, which is a polynomial 
in $n$ of degree $r$. The complexity bound follows from this and 
Corollary 2.4.
{\hfill \hfill \hfill} \qed
\enddemo

\head 4. Permanent Estimators \endhead 

Let us fix an ${\Bbb R}$-algebra $\AA$. We assume that as a real vector 
space, $\AA$ is finite-dimensional, so as a vector 
space $\AA$ is isomorphic to ${\Bbb R}^m$ for some $m$.
Suppose that there is a scalar 
product $\langle \cdot, \cdot \rangle: \AA \times \AA \longrightarrow 
{\Bbb R}$. In other words, for every two elements $a, b \in \AA$, a real 
number $\langle a, b \rangle$ is defined, such that 
\medskip
$\langle a, b \rangle = \langle b, a \rangle$;
\smallskip
$\langle \alpha a + \beta b, c \rangle =\alpha \langle a,c \rangle + 
\beta \langle b, c \rangle$ for $a,b,c \in \AA$ and 
$\alpha, \beta \in {\Bbb R}$ and
\smallskip
$\langle a, a \rangle \geq 0$ for all $a \in \AA$.  
\medskip
For example, if $\AA=\Mat(d, {\Bbb R})$ is a matrix algebra, one may 
choose $\langle a, b \rangle = \tr(ab^t)$.
Generally, we don't assume any relation between the multiplication in 
$\AA$ and the scalar product. We also allow the form
$\langle \cdot, \cdot \rangle$ to be degenerate,
that is, we allow $\langle a, a \rangle=0$ for some non-zero $a$.
As usual, we define $\|a\|^2=\langle a, a \rangle$.

Suppose further, that there is a Borel probability measure $\mu$ on $\AA$ 
with the properties:
\medskip
$$\EE x =\int_{\AA} x \ d\mu =0 \tag4.1$$ 
and for any (non-commutative) polynomial 
$$p(x_1, \ldots, x_n)=\sum_{1 \leq i_1, \ldots, i_n \leq n} 
\gamma_{i_1, \ldots, i_n} x_{i_1} \cdots x_{i_n}, \quad 
\text{where} \quad \gamma_{i_1, \ldots, i_n} \in {\Bbb R}$$
we have 
$$\EE \langle p(x_1, \ldots, x_n), p(x_1, \ldots, x_n) \rangle < +\infty,
\tag4.2$$
provided $x_1, \ldots, x_n \in \AA$ are chosen independently and 
at random.
\proclaim{(4.3) Theorem} Let $\AA$ be an ${\Bbb R}$-algebra with a 
scalar product $\langle \cdot, \cdot \rangle$ and a probability measure 
$\mu$ satisfying (4.1)--(4.2). Let us define a constant $c(n, \mu) \geq 0$ 
as follows:
 $$c(n, \mu)=\EE \|Z\|^2 , \quad \text{where} \quad 
Z={1 \over n!} \sum_{\sigma \in S_n} Y_{\sigma(1)} \cdots Y_{\sigma(n)}$$
and $Y_1, \ldots, Y_n$ are sampled independently and at random from 
$\AA$. 

For a given real non-negative matrix 
$A=(a_{ij})$, let us define a random $n \times n$ matrix $B \in \Mat(n, \AA)$,
$B=(b_{ij})$ as follows:
$$b_{ij} =\sqrt{a_{ij}} X_{ij},$$
where $X_{ij}$ are sampled independently and at random from $\AA$.
Then 
$$\EE \|\sdet B\|^2=c(n, \mu) \per A.$$
\endproclaim
\demo{Proof} By Definition 3.1,
$$\sdet B= {1 \over n!} \sum_{(\sigma, \tau) \in S_n \times S_n} 
(\sgn \sigma) (\sgn \tau) 
\biggl(\prod_{i=1}^n a_{\sigma(i) \tau(i)}\biggr)^{1/2}
X_{\sigma(1) \tau(1)} \cdots X_{\sigma(n) \tau(n)}.$$
Hence $\|\sdet(B)\|^2$ can be written as the sum of 
the terms
$$ \split &\Bigl({1 \over n!}\Bigr)^2 
(\sgn \sigma_1) (\sgn \tau_1) (\sgn \sigma_2) (\sgn \tau_2)
\biggl(\prod_{i=1}^n a_{\sigma_1(i) \tau_1(i)}\biggr)^{1/2}
\biggl(\prod_{i=1}^n a_{\sigma_2(i) \tau_2(i)}\biggr)^{1/2} \times \\ &
\big\langle X_{\sigma_1(1) \tau_1(1)} \cdots X_{\sigma_1(n) \tau_1(n)}, \ 
X_{\sigma_2(1) \tau_2(1)} \cdots X_{\sigma_2(n) \tau_2(n)} \big\rangle
\endsplit
\tag4.3.1$$
for all 4-tuples $(\sigma_1, \sigma_2, \tau_1, \tau_2) \in 
S_n \times S_n \times S_n \times S_n$.

Suppose that the sets of indices 
$$\Bigl\{\bigl(\sigma_1(i), \tau_1(i)\bigr):
i=1, \ldots, n \Bigr\} \quad \text{and} \quad
\Bigl\{\bigl(\sigma_2(i), \tau_2(i)\bigr): i=1, \ldots, n \Bigr\}$$ 
coincide. 
This is the case if and only if there is a  
permutation $\pi \in S_n$ such that $\sigma_2(i)=\sigma_1(\pi(i))$ and 
$\tau_2(i)=\tau_1(\pi(i))$ for $i=1, \ldots, n$, in which case 
the term (4.3.1) can be written 
as 
$$\split &\Bigl({1 \over n!}\Bigr)^2 \prod_{i=1}^n a_{\sigma_1(i) \tau_1(i)}
\times \\ &
\big\langle X_{\sigma_1(1) \tau_1(1)} \cdots X_{\sigma_1(n) \tau_1(n)}, \ 
X_{\sigma_1(\pi(1)) \tau_1(\pi(1))} \cdots 
X_{\sigma_1(\pi(n)) \tau_1(\pi(n))} \big\rangle. \endsplit$$
Let us fix $\sigma_1, \tau_1 \in S_n$ and let $\pi$ range over 
$S_n$. Then the expected value of the corresponding sum of monomials is
$${c(n, \mu) \over n!} \prod_{i=1}^n a_{\sigma_1(i) \tau_1(i)}.$$
If the sets of indices 
$$\Bigl\{\bigl(\sigma_1(i), \tau_1(i)\bigr):
i=1, \ldots, n \Bigr\} \quad \text{and} \quad
\Bigl\{\bigl(\sigma_2(i), \tau_2(i)\bigr): i=1, \ldots, n \Bigr\}$$ in 
(4.3.1) do not coincide, then there is a term $X_{ij}$, which is 
in the right hand side of the scalar product, but not in the left hand
side. The conditional expectation with respect to $X_{ij}$ is 0 because
of (4.1) and the linearity of expectation. Hence the expectation of 
the corresponding term is 0. Summarizing, the expected value of 
$\|\sdet B\|^2$ is 
$$\sum_{(\sigma_1, \tau_1) \in S_n \times S_n} 
{c(n, \mu) \over n!} \prod_{i=1}^n a_{\sigma_1(i) \tau_1(i)}=
c(n, \mu) \per A.$$
{\hfill \hfill \hfill} \qed
\enddemo
Algorithm 3.5 and Theorem 3.6 ensure the property (1.2.1) of the estimator.
Theorem 4.3 implies (1.2.3) and (1.2.2) (up to a constant $c(n, \mu)$;
for sufficiently generic measure and scalar product
the constant is strictly positive). 
One can see that (1.2.4) is satisfied as well.

\head 5. A Series of Estimators Conjectured to be Asymptotically 
Exact \endhead

Let us fix a positive integer $d$ and let $\AA=\Mat(d, {\Bbb R})$ be 
the algebra of real $d \times d$ matrices. We introduce the scalar 
product $\langle \cdot, \cdot \rangle $ on $\AA$ by letting
$$\langle a, b \rangle = \tr(a b^t).$$
Let $\mu$ be the standard Gaussian measure on $\AA$ with the density
$$\psi(a)=(2 \pi)^{-d^2/2} e^{-\|a\|^2/2}.$$
Thus to sample a random matrix $a \in \AA$, we  sample each of its 
$d^2$ entries independently from the standard normal distribution.

The construction of Section 4 gives us a permanent estimator.
Although we don't know the value of $c(n, \mu)$, we don't really need
it, since we can approximate it by applying the same algorithm to 
the identity matrix.

Hence, for each positive integer $d$, we obtain a permanent estimator. 
\subhead (5.1) The $d$-th estimator \endsubhead 
\bigskip
\noindent{\bf Input:} A non-negative $n \times n$ 
real matrix $A=(a_{ij})$.
\medskip
\noindent{\bf Output:} A non-negative number $\alpha$ approximating
$\per A$.
\medskip
\noindent{\bf Algorithm:} Sample $n^2$ matrices $u_{ij}$ of size $d \times d$ 
 from the standard Gaussian distribution in $\Mat(d, {\Bbb R})$.
Define $n \times n$ matrices $B=(b_{ij})$, $E=(e_{ij})$ with 
$b_{ij}, e_{ij} \in \Mat(d, {\Bbb R})$ for $i,j=1, \ldots, n$ 
as follows:
$$b_{ij}=u_{ij}\sqrt{a_{ij}} \quad \text{and} \quad e_{ij}=\cases u_{ii} 
&\text{if\ } i=j \\ 0 &\text{if\ } i \ne j. \endcases$$
Apply Algorithm 3.5 to compute $d \times d$ matrices $\sdet B$ and 
$\sdet E$. 
Compute 
$$ \alpha = \|\sdet B \|^2 / \|\sdet E\|^2 .$$
Output $\alpha$.  
\bigskip
We conjecture that the output $\alpha$ approximates $\per A$ within 
an exponential factor $\gamma_d^n$ and that 
$\gamma_d$ approaches 1 as $d$ grows. More precisely, we conjecture that
there exist 
\medskip
a sequence $\{\gamma_d\}$ of non-negative real numbers such that 
$\lim_{d \longrightarrow +\infty} \gamma_d=1$
\smallskip
and 
\smallskip
a sequence of functions $\{f_d(n, \epsilon)\}$ 
such that for any $\epsilon>0$ and any $d$ 
we have 
$\lim_{n \longrightarrow +\infty} f_d(n, \epsilon)=0$;
\medskip
\noindent so that for any $n \times n$ non-negative matrix $A$ and any 
$\epsilon>0$, the 
output $\alpha$ of the $d$-th estimator (5.1) satisfies:
$$\Pr \Bigl\{ (\gamma_d -\epsilon)^n \per A \leq  \alpha
\leq  (\gamma_d +\epsilon)^n 
\per A \Bigr\} \geq 1 -f_d(n, \epsilon).$$
Note, that the complexity of the $d$-th algorithm is $O(n^{d^2+3})$.
We discuss some plausible reasons why the conjecture might be true.
\subhead (5.2) Why there might be sharp concentration \endsubhead

From Theorem 4.3 and Chebyshev's inequality, one deduces that 
the value of 
$c^{-1}(n, \mu)\|\sdet B\|^2$ is unlikely to overestimate 
$\per A$:
$$\Pr \bigl\{ c^{-1}(n, \mu)\|\sdet B\|^2 \geq K \per A 
\bigr\} \leq 1/K$$
for every $K>0$.
Hence the main problem is to prove that 
$c^{-1}(n, \mu)\| \sdet B\|^2$ is unlikely to underestimate
$\per A$ as well.

Suppose that in Algorithm 5.1, instead of the symmetrized determinant of $B$,
we take the Cayley determinant $\Cdet$ (see Section 1.3). One can prove
that 
$$\EE  \|\Cdet B\|^2 = d^n \per A.$$
Moreover, the method of [Barvinok 99] carries over and one can prove that for 
any $1>\epsilon>0$
$$\Pr \Bigl\{d^{-n} \|\Cdet B\|^2 \leq 
(\epsilon {\frak c}_d)^n  \per A \Bigr\} \leq {8 \over n \ln^2 \epsilon},$$
where constant ${\frak c}_d$ is defined as follows (see [Barvinok 99]):
let $\xi_1, \ldots, \xi_d$ be independent random variables having 
the standard Gaussian density $(2 \pi)^{-1/2} e^{-x^2/2}$. Then 
$${\frak c}_d = \exp\biggl\{ \EE \ln \Bigl({\xi_1^2 + \ldots + \xi_d^2 \over 
d} \Bigr) \biggr\}.$$  
We have $\lim_{d \longrightarrow +\infty} {\frak c}_d=1$ and, in fact,
${\frak c}_d=1+O(1/d)$.

The symmetrized determinant $\sdet B$ can be considered as the 
average of the $n!$ Cayley determinants of the matrices obtained from 
$B$ by permuting rows in all possible ways. It seems quite plausible 
to the author that the concentration for $\|\sdet B\|^2$
should be at least as sharp as for $\|\Cdet B\|^2$.  

\head Acknowledgment \endhead

I am grateful to Rishi Raj for many helpful discussions.

\head References \endhead   

\noindent [Aslaksen 96] H. Aslaksen,  Quaternionic determinants,
{\it The Mathematical}

{\it Intelligencer}, {\bf 18}(1996), 57--65.
\smallskip
\noindent [Bapat and Raghavan 97] R.B. Bapat and T.E.S. Raghavan,
{\it Nonnegative Matrices}

{\it and Applications}, Encyclopedia of Mathematics
and its Applications, {\bf 64}, 

Cambridge Univ. Press, 
Cambridge, 1997.
\smallskip
\noindent [Barvinok 97] A. Barvinok, 
Computing mixed discriminants, mixed volumes, and 

permanents,
{\it Discrete $\&$ Computational Geometry}, {\bf 18} (1997), 205--237.
\smallskip
\noindent [Barvinok 99] 
A. Barvinok, Polynomial time algorithms to approximate permanents 

and mixed
discriminants within a simply exponential factor, 

{\it Random Structures  $\&$ Algorithms}, {\bf 14}(1999), no. 1, 29--61. 
\smallskip
\noindent [Gelfand and Retakh 97],
I. Gelfand and V. Retakh, Quasideterminants. I. 

{\it Selecta Mathematica (N.S.)}, {\bf 3}(1997), 517--546. 
\smallskip
\noindent [Girko 90] V.L. Girko, {\it Theory of Random Determinants},
Mathematics and 

its Applications, {\bf 45}, Kluwer, Dordrecht, 1990. 
\smallskip
\noindent [Jerrum and Sinclair 89] 
M. Jerrum and A. Sinclair, Approximating the permanent, 

{\it SIAM Journal on Computing}, {\bf 18} (1989), 
1149--1178.
\smallskip
\noindent [Karmarkar {\it et al.} 93]
N. Karmarkar, R. Karp, R. Lipton, L. Lov\'asz and M. Luby,

A Monte Carlo algorithm for estimating the permanent, {\it SIAM Journal}

{\it on Computing}, {\bf 22}(1993), 284--293.
\smallskip
\noindent [Linial {\it et al.} 98] 
N. Linial, A. Samorodnitsky and A. Wigderson, A deterministic 

strongly 
polynomial algorithm for matrix scaling and approximate permanents,

{\it Proc. 30 ACM Symp. on Theory of Computing}, ACM, New York, 1998, 644--652;

revised version is in {\it Combinatorica}, to appear.
\smallskip
\noindent [Lov\'asz and Plummer 86] L. Lov\'asz and M.D. Plummer, 
{\it Matching Theory}, 

North - Holland, Amsterdam - New York and
 Akad\'emiai Kiad\'o, Budapest, 1986.
\smallskip
\noindent [Papadimitriou 94] C.H. Papadimitriou, {\it Computational 
Complexity},

Addison-Wesley, Reading, Mass., 1994.
\smallskip
\enddocument

\end